\documentclass[10pt]{article}
\usepackage{amsmath,amsthm,hyperref}
\usepackage[latin1]{inputenc}
\usepackage{picinpar}
\usepackage{longtable}
\usepackage{amsrefs}
\usepackage[dvips]{graphicx}
\usepackage[mathscr]{eucal}
\usepackage{calc}
\usepackage{ifthen}
\usepackage{dcpic,pictexwd}
\usepackage[all]{xy}
\usepackage{enumerate}
\usepackage{amssymb,latexsym}
\usepackage{amsthm}
\usepackage[dvips]{graphics}
\usepackage{color}
\usepackage{tabularx}
\usepackage{makeidx}
\usepackage{tkz-graph}
\usepackage{tikz-cd}
\usetikzlibrary{calc}
\usepackage{tikz}
\usetikzlibrary{arrows}
\usepackage{stackrel}
\usepackage{multicol}
\usepackage{float}
\newtheorem{teor}{Theorem}
\newtheorem{lema}{Lemma}
\newtheorem{prop}{Proposition}
\newtheorem{corol}{Corollary}
\newtheorem{defin}{Definition}

\newtheorem{Exp}{Example}
\newtheorem{rem}{Remark}

\DeclareMathOperator{\rep}{Rep}
\DeclareMathOperator{\repq}{Repq}

\begin{document}
\thispagestyle{plain}
\par\bigskip
\begin{centering}

\textbf{On the Auslander-Reiten quiver for the category of representations of partially ordered sets with an involution.}

\end{centering}\par\bigskip
\begin{centering}
\footnotesize{Raymundo Bautista Ramos}\\
\footnotesize{Verónica Cifuentes Vargas}\\
\end{centering}
\par\bigskip
\small{In this article we describe the Auslander-Reiten quiver for some posets with an involution, that we call types $\mathfrak{U}_n$ and $\mathfrak{U}_\infty$. These posets appear in the differentiation III of Zavadskij [see, \cite{Zavadskij12}]. We follow the approach to classical Auslander-Reiten theory due to Auslander, Reiten and Smal\o{} \cite{Auslander}. For this purpose, we give a natural exact structure for the category of representations of a partially ordered set with an involution. We describe the
	projective and injective representations. }


\section{Introduction}\label{sec:1}
Matrix representations of a partially ordered set (poset) with an equivalence relation were introduced by Nazarova and Roiter in \cite{Nazarova2}. A particular but important case is when each equivalence class consits of at most two elements. In this case we say that the poset has an involution. Following the ideas presented by Gabriel in \cite{Gabriel1}, Zavadskij (see section 9 of \cite{Zavadskij12}) consider filtered $k$-linear representations of a poset with an involution ($k$ a field) with the purpose of  giving a  strict foundation to the several differentiation algorithms introduced by him. These algorithms were conceived in terms of matrices. \par\bigskip

The different differentiation algorithms are a powerful tool for classifing posets with an involution of finite and tame type.  These algorithms  can be applied when there is an appropiated subposet  of the poset with an involution, then the differentiation consists in the change of this subposet obtaining a new poset.  After this, is defined a  functor from the representations  of the original poset  into the new one in such a way that with the exception of a finite number of isomorphism clases of indecomposables, one obtains a bijection between the isomorphism clases of  the indecomposable representations of  the original poset and the isomorphism clases of indecomposables of the new poset.

\par\bigskip
For some cases it is well known that the functor given by the differentiation,  induces an equivalence between the corresponding categories of representations modulo the ideal generated by a finite number of indecomposable objects.  Here we are interested in the appropiated subposet appearing  in the differentiation III of Zavadskij. We will describe its Auslander-Reiten quiver. To do that we first make some general considerations on a natural exact structure on the category of vector space representations of a partial ordered set with an involution  and their Auslander-Reiten sequences.
\par\bigskip
In section 3 we will see that the above exact structure has enough projectives. Taking
the endomorphism ring of the direct sum of a representative set of the indecomposable
projectives we obtain a right-peak algebra in the sense of Simson (see, \cite{Simson1}). The category of representations of the poset with involution is equivalent with
the category of the socle projective modules of this last algebra. Posets with an equivalence 
relation are a special case of completed posets introduced by Nazarova and Roiter. A further generalization is the notion of stratified posets introduced by Simson (see, 17.8 of \cite{Simson}). A vector space category in the sense of Ringel (see, \cite{Ringel}) is associated to a stratified poset. Given a vector space category $\mathbb{K}_F$, there are two equivalent categories of representations,
the subspace category $\mathcal{U}(\mathbb{K}_F )$, and $\mathcal{V}(\mathbb{K}_F )$ the factor space category. Moreover a right-peak algebra $\mathbb{R}_{\mathbb{K}_F}$
is associated to any vector space category $\mathbb{K}_F$. The category $\text{mod}_{sp}(\mathbb{R}_{\mathbb{K}_F})$ is
equivalent to the category $\mathcal{V}(\mathbb{K}_F )$ modulo those morphisms which are factorized through
finite direct sums of a finite familiy of objects. In our case if $\mathbb{K}_{F}$ is the vector space category associated to a poset with involution, the endomorphism algebra of the projective generator of the category of representations of our poset with involution is isomorphic to
$\mathbb{R}_{\mathbb{K}_F}$.
\par\bigskip
Now, we give the layout of the contents of this paper.
In section 2 we introduce the main
definitions. In section 3 we present an exact structure for the category of representations of
a partially ordered set with involution. In the second part of this section we consider representations by factor spaces instead of subspaces. This allows us to describe the injective
objects in the category of representations. In section 4 we take a projective generator for
the category of representations of a poset with involution, obtaining a right-peak algebra.
We prove that the category of the socle-projective modules of this algebra is equivalent
to the category of representations of the poset with involution. Finally in section 5 we
describe the Auslander-Reiten quiver of a poset with involution $\mathfrak{U}_n$ and its extension $\mathfrak{U}_\infty$.
\section{Partially ordered set with an involution}
In this section we define partially ordered sets with an involution and their category of representations.

\begin{defin}
	A partially ordered set with an equivalence relation is a triple $(\mathscr{P}, \leq, \theta)$, where  $(\mathscr{P}, \leq)$ is a partially ordered set and in  $\mathscr{P}$ there is an equivalence relation whose equivalence classes is $\theta$. If the cardinality of each equivalence class is less than or equal to two, we will say that triple  $(\mathscr{P}, \leq, \theta)$ is a partially ordered set with an involution.  If $x\in \mathscr{P}$ we will denote by $[x]$ its equivalence class.
\end{defin}

\begin{rem}
	
	From now on we will omit the order relation in the notation for poset with an involution, that is, we will write $(\mathscr{P},\theta)$ instead of  $(\mathscr{P}, \leq, \theta)$.
\end{rem}

\begin{Exp}
	Let $(\mathscr{P},\theta)$ be a poset with involution  where  $\mathscr{P}$ is as in Figure 1 with $a <b$, $c<b$, $c<a^{*}$ and $\theta=\{(a,a^{*}),b,c\}$.
	\begin{figure}[H]
		\begin{center}
			\begin{tikzpicture}
				
				\node (1) at (-1.3,-0.4) {$(\mathscr{P},\theta)$=};
				\node (1) at (-0.3,0.11) {$b$};
				\node (1) at (-0.3,-1) {$a$};
				\node (1) at (1.3,0.13) {$a*$};
				\node (1) at (1.3,-1) {$c$};
				
				\node (1) at (0,0) {$\circ$};
				\node [right  of=1, node distance=1cm] (2)  {$\bullet$};
				
				\node [below  of=1, node distance=1cm] (1')  {$\bullet$};
				\node [right  of=1',node distance=1cm] (2')  {$\circ$};

				\draw [ shorten <=-2.2pt, shorten >=-2.2pt] (1') -- (1);
				\draw [  shorten <=-2pt, shorten >=-2pt] (2') -- (2);
				\draw [  shorten <=-2pt, shorten >=-2pt] (2') -- (1);
			\end{tikzpicture}	
		\end{center}
		\caption{Diagram of a poset with an involution}
	\end{figure}
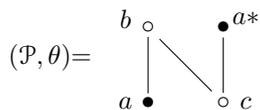
\end{Exp}

\subsection{Vector space representations for posets with an  involution}{\label{vecRep}}

Zavadskij  introduced filtered $k$-linear representation of posets with an  involution $(\mathscr{P},\theta)$  \cite{Zavadskij12}. Here we introduce an equivalent  definition to the one given by him. For this, we consider $(\mathscr{P},\theta)$  a poset with an involution.  We take $V_0$ a $k$-vector space and $z \in \theta$, take $V_0^{z}$ the $k$-vector space consisting   of all functions $h: z \xrightarrow{} V_0$. For $x \in z$, we have the inclusion:
$i_x: V_0 \xrightarrow{} V_0^{z}$, defined by	$$ i_x(v)(y)=\left\{ \begin{array}{lcc}
	0, & \hspace{-0.3cm}  \textnormal{if}  & y \neq x, \\
	\\ v, & \quad  \quad \quad \textnormal{otherwise}. & \\
\end{array}
\right. $$
and the projection in the summand $x$ of $V^{z}_0$,
$\pi_x:V_{0}^{z}\to V_0$,  that is, for $h\in V^{z}_0$, $\pi_x(h)=h(x)$.\\

In the following,  if $V$ is a $k$-vector subspace of $ V_0^{z}$
and $x\in z$, 
\begin{align*} 
	V_{x}^{-} &= i^{-1}_x(V)= \{v\in V_0\mid i_x(v)\in V\},\\
	V_x^+ &= \pi_x(V)=\{h(x)\mid h\in V\}.
\end{align*}

\smallskip
\begin{defin}{\label{defPI}}
	A vector space representation $V=(V_0,V_z)_{z\in\theta}$ of $(\mathscr{P, \theta})$ is given by:
	\begin{enumerate}
		\item a finite-dimensional $k$-vector space $V_0$,
		\item for each $z\in \theta$, a vector subspace $V_z$ of $V^{z}_0$ such that if  $y<x$ then 
		$$V^{+}_y\subset V^{-}_x.$$
	\end{enumerate}
\end{defin}

\begin{Exp}{\label{RepV}}
	Let $(\mathscr{P},\theta)$ be a poset with an  involution  where  $\mathscr{P}$ is as in  Figure 2 with $a <b^{*}$, $a^{*}<b$, $a^{*}<b^{*}$ and $\theta=\{(a,a^{*}),(b,b^{*})\}$.\smallskip
	
	\begin{figure}[H]
		\begin{center}
			\begin{tikzpicture}
				
				\node (1) at (-1.3,-0.4) {$(\mathscr{P},\theta)$=};
				\node (1) at (-0.3,0.11) {$b^*$};
				\node (1) at (-0.3,-1) {$a$};
				\node (1) at (1.3,0.13) {$b$};
				\node (1) at (1.3,-1) {$a^*$};
				
				\node (1) at (0,0) {$\bullet$};
				\node [right  of=1, node distance=1cm] (2)  {$\bullet$};
				
				\node [below  of=1, node distance=1cm] (1')  {$\bullet$};
				\node [right  of=1',node distance=1cm] (2')  {$\bullet$};

				\draw [ shorten <=-2.2pt, shorten >=-2.2pt] (1') -- (1);
				\draw [  shorten <=-2pt, shorten >=-2pt] (2') -- (2);
				\draw [  shorten <=-2pt, shorten >=-2pt] (2') -- (1);
				\draw [ shorten <=-2pt, shorten >=-2pt] (2) -- (1');
			\end{tikzpicture}
		\end{center}
		\caption{Diagram of a poset with an  involution.}
	\end{figure}
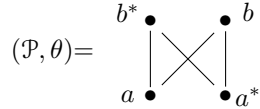
	
	We will show that  $V=(V_0, V_{(a,a^*)}, V_{(b,b^*)})$ is a vector space representation  of   $(\mathscr{P},\theta)$,   where $V_0=\mathbb{R}^3$, $\mathscr{B}=\{e_1,e_2,e_3\}$ is the canonical basis of $V_0$ and  $V_{(a,a^*)}=\langle h\rangle$, with
	
	\begin{center}
		$\begin{array}{cccc}
			h:& (a,a^*) &\rightarrow \quad  \mathbb{R}^3 & \\
			& a   & \hspace*{-1cm}\mapsto & \hspace*{-1.3cm}e_1\\
			& a^*   & \hspace*{-1cm}\mapsto & \hspace*{-1.3cm}e_2\\
		\end{array}
		$
	\end{center}
	
	and $V_{(b,b^*)}=\langle h_1,h_2,h_3,h_4\rangle$, with
	
	\begin{center}
		$\begin{array}{cccc}
			h_1:& (b,b^*)&\rightarrow \quad  \mathbb{R}^3 & \\
			& b   & \hspace*{-1cm}\mapsto & \hspace*{-1.3cm}e_2\\
			& b^*   & \hspace*{-1cm}\mapsto & \hspace*{-1.3cm} 0\\
		\end{array}
		$ \hspace{1cm}
		$\begin{array}{cccc}
			h_2:& (b,b^*) &\rightarrow \quad  \mathbb{R}^3 & \\
			& b   & \hspace*{-1cm}\mapsto & \hspace*{-1.3cm}0\\
			& b^*   & \hspace*{-1cm}\mapsto & \hspace*{-1.3cm}e_1\\
		\end{array}
		$\\
		
		\vspace{0.8cm}
		
		$\begin{array}{cccc}
			h_3:& (b,b^*) &\rightarrow \quad  \mathbb{R}^3 & \\
			& b   & \hspace*{-1cm}\mapsto & \hspace*{-1.3cm}0\\
			& b^*   & \hspace*{-1cm}\mapsto & \hspace*{-1.3cm}e_2\\
		\end{array}$  
		\hspace{1.cm}
		$\begin{array}{cccc}
			h_4:& (b,b^*) &\rightarrow \quad  \mathbb{R}^3 & \\
			& b   & \hspace*{-1cm}\mapsto & \hspace*{-1.3cm}e_3\\
			& b^*   & \hspace*{-1cm}\mapsto & \hspace*{-1.3cm}e_3\\
		\end{array}
		$
	\end{center}
	
	Indeed, 	$V_a^{+}=\pi_a(V_{(a,a^*)})=\mathbb{R}\{e_1\}$,    $V_{a^*}^{+}= \pi_{a^*}(V_{(a,a^*)})=\mathbb{R}\{e_2\}$,
	$V_b^{-}=i_{{b}}^{-1}(V_{(b,b^*)})=\mathbb{R}\{e_2\}$  and  
	$\quad  V_{b^*}^{-}=i_{{b}^{*}}^{-1}(V_{(b,b^*)})=\mathbb{R}\{e_1,e_2\}$;
	then, for  $a<b^*$ is obtained that  $V_a^{+}\subset V_{b^*}^{-}$. For $a^*<b^*$,   it is satisfied that   $V_{a^*}^{+}\subset V_{b^*}^{-}$ and  for  	$a^*<b$, it is true that
	$V_{a^*}^{+}\subset V_{b}^{-}$.
\end{Exp}
\smallskip

\begin{defin}{\label{homo}}
	If $ V=(V_0, V_z)_{z\in  \theta}$ and $W=(W_0,W_z)_{z\in \theta}$ are two representations
	in  $(\mathscr{P}, \theta)$, and $\varphi:V_0\to W_0$ is a morphism of vector spaces, such that for each $z\in \theta$, we have the morphism $\varphi^z:V_0^z\to W_0^z$ and  for $h:z\to V_0$, $ \varphi^z(h)=\varphi h$. Then  a morphism  $ V\to W$ consists of a morphism of vector space  $\varphi:V_0\to W_0$ such that,  for all $z\in \theta$
	\begin{center}
		$\varphi^z(V_z) \subset W_z.$
	\end{center}
\end{defin}

\smallskip

With the definitions above, the vector space representations of a poset with an involution can also be viewed categorically. We denote by  $\textnormal{Rep}(\mathscr{P},\theta)$ this category.

\begin{defin}{\label{directa}}
	If $ V=(V_0,V_z)_{z\in  \theta}$ and $W=(W_0,W_z)_{z\in \theta}$ are two representations
	of $\textnormal{Rep}(\mathscr{P}, \theta)$ then their direct sum is 
	$$V\bigoplus W= (V_0\bigoplus W_0, V_z \bigoplus W_z )_{z\in \theta}.$$
\end{defin}

\section{Exact structure on the category  $\rep(\mathscr{P},\theta).$}{\label{sec2.1}}

Let  $\mathscr{A}$ be an additive category in which all idempotents split, and let $\varepsilon$  be  
a collection of pairs of morphisms  $ M \xrightarrow{u} E \xrightarrow{v} N$. A morphism $u:M\xrightarrow{u} E$ is called an $\varepsilon$-\textit{inflation} if there exists a morphism  $v:E\to N$ such that $(u,v)\in \varepsilon$. A morphim $v:E\xrightarrow{v} N$ is called an  $\varepsilon$-\textit{deflation} if there exists a morphism  $u: M  \to E$ such that $(u,v)\in \varepsilon$.

The pair $(\mathscr{A}; \varepsilon)$ will be called \textit{exact structure} (see, \cite{Gabriel2}, \cite{Reiten}, and  \cite{Buhler}) if the following conditions are satisfied:
\begin{enumerate}
	\item  The family $\varepsilon$ is closed under 
	isomorphisms; that is, if there exists a commutative diagram:
	$$
	\begin{tikzcd}
		M  \arrow[r, "u"]\arrow[d, "s"] & E \arrow[r, "v"]\arrow[d, "t"]
		& N \arrow[d, "r"] \\
		M_1  \arrow[r, "u_1"] & E_1 \arrow[r, "v_1"]
		& N_1 
	\end{tikzcd}$$
	
	where $s, t, r$ are  isomorphisms and  the top row is in $\varepsilon$, then the bottom  row
	belongs to $\varepsilon$.
	\item  If $(u,v)\in \varepsilon$, then  $u$ is a kernel of $v$  and $v$  is  a  cokernel of $u$.
	\item $\text{id}_M: M \to M$  is both  $\varepsilon$-inflation
	and  $\varepsilon$- deflation.
	\item  	\begin{description}
		\item[a.]For each $\varepsilon$-sequence $M \xrightarrow{f} E\xrightarrow{g} N $ and each morphism $w: X \to N $ there are morphisms $\beta: F \to X$  and $ \lambda: F \to E$ such that the pair $(\lambda,\beta )$ is a pullback of the pair $(g,w)$ and $\beta$ is an $\varepsilon$-deflation.
		\item[b.]  For each $\varepsilon$-sequence $M \xrightarrow{f} E\xrightarrow{g} N $ and each morphism $u: M \to X $ there are morphisms $\alpha: X \to F$  and $ \lambda: E \to F$ such that the pair $(\alpha, \lambda)$ is a pushout of the pair $(u,f)$ and $\alpha$ is an $\varepsilon$-inflation.
	\end{description}
	\item The composition of $\varepsilon$-inflations ($\varepsilon$-deflations, respectively) is again an  $\varepsilon$-inflation ($\varepsilon$-deflation, respectively).
	\item If $u_2u_1$ is an $\varepsilon$-inflation then $u_1$ is an $\varepsilon$-inflation. If $v_2v_1$ is an $\varepsilon$-deflation then $v_2$ is an $\varepsilon$-deflation.
\end{enumerate}

In our case, let $(\mathscr{P},\theta)$ be a poset with an  involution and let us consider $\varepsilon$ the family of sequences of morphisms:
$$(V_0, V_z) \xrightarrow{u} (E_0, E_z) \xrightarrow{v} (W_0,W_z),$$
in the category $\rep(\mathscr{P},\theta)$ such that:
\begin{enumerate}
	\item The sequence $0  \to  V_0  \xrightarrow{u} E_0 \xrightarrow{v} W_0\to 0 $ is exact.
	\item For all $z\in \theta$, the sequence  $0  \to  V_z \xrightarrow{u_z} E_z \xrightarrow{v_z} W_z\to 0 $ is exact.
\end{enumerate}



From now on we will denote a representation $(U_0, U_z)_{z\in \theta}$ in $\rep(\mathscr{P},\theta)$ by
$U$, unless it is necessary to specify the subspaces $U_z$.  
\begin{defin}{\label{epim}} 
	A morphism $f: V\to  W$ in the category $\rep(\mathscr{P}, \theta)$ will be called a
	proper epimorphism if $f: V_0\to  W_0$ is an epimorphism and for each $z\in \theta$,
	$f^z: V_{0}^z\to W_{0}^z$ induces an epimorphism $f^z: V_z\to W_z$.
\end{defin}

\begin{prop}{\label{sequen1}}
	Let  $f: V\to  W$ be a proper epimorphism. If  $U_0 = \textnormal{ker}(f)$  and  $U_z = \textnormal{ker}(f^z) \cap V_z$ then $U=(U_0, U_z)_{z\in \theta}$ is a representation of $(\mathscr{P},\theta)$. 
\end{prop} 
\textbf{Proof}. In the first place, we observe that for each $z\in \theta$,  $\text{ker}(f^z) =
U_{0}^z$, then $U_z = U_{0}^z \cap V_z$. We suppose now that $a <a_1$  and $ (x, y) \in  U_{(a,b)}$ so,  $(x, y) \in  V_{(a,b)}$  then if  $z_1 = (a_1, b_1), (x, 0) \in  V_{(a_1,b_1)}$  as, $(x, y) \in  U_{(a,b)}$,
$(f(x), f(y)) = (0, 0)$ then  $f(x) = 0$; therefore $(x, 0) \in  U_{(a_1,b_1)}$. This proves
that in effect $(U_0, U_z)_{z\in \theta}$ is a representation of  $(\mathscr{P},\theta)$.
$\hfill\blacksquare$
\smallskip
\begin{corol}
	If $f: V\to W$ is a proper epimorphism and 
	$U$ is as in  the  previous proposition,  then  an $\varepsilon$-sequence 
	$$U \xrightarrow{g} V\xrightarrow{f} W,$$
	
	is obtained. Therefore, $f$ is an $\varepsilon$-deflation if and only if $f$ is a proper epimorphism.
\end{corol}

\begin{defin}{\label{monom}}
	A morphism $f: U \to V$ in the category $\rep(\mathscr{P}, \theta)$ will be called a proper
	monomorphism if $f: U_0 \to V_0$ is a monomorphism and for each $z \in \theta$, $f^z(U_z)= U_0 \cap f(U)$.
\end{defin}

\begin{prop}{\label{sequen2}}
	
	If $f: U \to V$  is a proper monomorphism then there exists a sequence
	in $\varepsilon$:
	$$U \xrightarrow{f} V \xrightarrow{g} W. $$
	Therefore $f$ is an $\varepsilon$-inflation  if and only if $f$ is a
	proper monomorphism.
\end{prop}
\textbf{Proof.}
Let  $g: V_0 \to W_0$ be  the cokernel of  $f$. For $z \in \theta$  we define $W_z = f^z(V_z)$. We will check that  $(W_0,W_z)_{z\in \theta}$ is a representation of $(\mathscr{P},\theta)$. Indeed,
let $(x, y) \in  W_{(a,b)}$ and  $(a_1, b_1) \in \theta$  with $a <a_1$.  Then $x = g(x_1)$, $y= g(y_1)$ with $(x_1, y_1)\in  V_{(a,b)}$ so,  $(x, 0)\in  V_{(a_1,b_1)}$, therefore $(x, 0) =
g(x_1, 0) \in  W_{(a_1,b_1)}$. This proves that  $(W_0,W_z)_{z\in \theta}$ is a representation. We prove  now that for each $z\in \theta$, the sequence:
$$0 \to U_z\xrightarrow{f^z} V_z \xrightarrow{g^z} W_z \to  0, $$
is exact. Since  $f^z$ is  a monomorphism, $g^z$ is an epimorphism and  $g^zf^z = 0$. It only remains to prove that if  $(x, y) \in  V_z$ is such that  $(g(x), g(y)) = (0, 0)$ then $(x, y) \in  f^z(U_z)$. Since the sequence
$$0 \to U \xrightarrow{f} V \xrightarrow{g} W \to 0, $$
is exact, then $(x,y)\in  f(U)$. Also $(x, y) \in  V_z$  and $f$ is a proper monomorphism, it  follows that  $(x, y) \in  f^z(U_z)$. This proves our claim.
$\hfill\blacksquare$
\par\bigskip

\smallskip
\begin{prop}{\label{Exac}}
	The pair $(\rep(\mathscr{P}, \theta), \varepsilon)$ is an exact category.
\end{prop}

\textbf{Proof.} 
Conditions  $1, 2$ and  $3$ are verified directly;   our characterization of $\varepsilon$-deflations and $\varepsilon$-inflations implies items 5 and 6. Thus, it remains to prove item 4. Let us show condition 4.\textbf{a.} Let 
$$U\xrightarrow{u} E \xrightarrow{v} V,$$	
be an $\varepsilon$-sequence and let
$f: W\to V$  be a morphism. Consider the morphism:
$$\phi=(v,-f):E\bigoplus W \to V,$$ here $v:W_0\to V_0$ and  $v^z:W_z\to V_z$ are epimorphisms for all $z\in \theta$. Then  $(v,-f):E_0\bigoplus W_0 \to V_0$ and  $\phi_z=(v^z,f^z):E_z\bigoplus W_z \to V_z$ 
are  epimorphisms for all $z\in \theta$, so  $\phi=(v,-f)$ is  a deflation, therefore there is an 
$\varepsilon$-sequence.
$$L\xrightarrow{(h1,-h2)^
	t} E\oplus W \xrightarrow{\phi} V.$$

This implies we have pull-back diagrams:
\[
\begin{tikzcd}
	L\arrow[r, "h_2"]\arrow[d, "h_1"] &  W \arrow[d, "f"]  \\
	E \arrow[r, "v"] & V 
\end{tikzcd}, \qquad  \begin{tikzcd}
	L_0\arrow[r, "h_2"]\arrow[d, "h_1"] &  W_0 \arrow[d, "f"]  \\
	E_0 \arrow[r, "v"] & V_0
\end{tikzcd}, \qquad \begin{tikzcd}
	L_z\arrow[r, "h^z_{2}"]\arrow[d, "h^z_{1}"] &  W_z \arrow[d, "f^z"]  \\
	E_z \arrow[r, "v^z"] & V_z
\end{tikzcd}
\]
\smallskip
here $v: E_0 \to V_0$ and $v^z:E_z \to V_z$ are epimorphisms, then $h_2: L_0\to  W_0$ and $h^z_{2}: L_z \to W_z$ are epimorphisms for all $z \in \theta$, this implies that $h_2$ is a deflation. This proves 4.\textbf{a}; item 4.\textbf{b}
is proved in a similar way.

$\hfill\blacksquare$

\subsection{ $\varepsilon$-projective representations}

\begin{defin}{\label{projec}}
	A representation $P$ of $(\mathscr{P},\theta)$  is called $\varepsilon$-projective if given an $\varepsilon$-deflation $g:E\to V$ and a morphism $f:P\to V$, there exists  a morphism $h:P\to E$ such that $gh=f$.
\end{defin}
\smallskip
\begin{rem}{\label{Rem2}}
	The representation ${S}=(S_0, S_z)_{z\in \theta}$
	with $S_0=k$ and $S_z=0$ for all $z\in \theta$, is a  projective representation. 	
\end{rem}
\smallskip
Let $w=(a,b)\in \theta$, we will define the representation $P(w)=(P(w)_0,P(w)_z)_{z\in \theta}$,  where $P(w)_0=k\langle e_1,e_2\rangle$ the  vector space of dimension two with bases $e_1,e_2$. If $a$ and $b$ are incomparable $P(w)_w=\langle (e_1,e_2)\rangle$, while if $a<b$ then $P(w)_w=\langle(0,e_1), (e_1,e_2)\rangle$.
\par\bigskip
Henceforth, we will use the following notation,  if $d_1,d_2 \in \mathscr{P}$ then
\begin{align*}
	\lambda(d_1,d_2)&=\begin{cases}
		1 &  \text{if} \quad  d_1<d_2,\\
		0 & \text{otherwise}.
	\end{cases}
\end{align*}

If $z = (a_1, b_1)$, $P(w)_z= \langle(\lambda(a,a_1)e_1,0), (0, \lambda(a, b_1)e_1), (\lambda(b, a_1)e_2, 0), (0, \lambda(b, b_1)e_2)\rangle.$\\

In case that, $w=\{a\}$ then $P(w)_0=k \langle e\rangle$ and for $z=(a_1,b_1)$ 
$$ P(w)_z =\langle (\lambda(a,a_1)e,0), (0, \lambda(a, b_1)e) \rangle.$$

It can be verified that  ${P}(w)=(P(w)_0, P(w)_z)_{z\in \theta}$   is in fact a representation.
\smallskip
\begin{defin}
	The element $(e_1, e_2) \in  P(w)$  will be  called the generator of the representation $(P(w)_0,P(w)_z)_{z\in \theta}$, when $w = (a, b)$ while the element $e \in  P(w)$  is the generator when  $w$  consists of a single point. 
\end{defin}

\smallskip
\begin{prop}{\label{cubierta}}
	Let $(V_0,V_z)_{z\in \theta}$ be a representation of  $(\mathscr{P},\theta)$, then if $w=(a,b)\in \theta$
	and $v\in V_w$ there exists an unique morphism  $f: P(w) \to  (V_0, V_z)_{z\in \theta}$ such that  $f((e_1, e_2)) = v$. If  $w =\{a\}$ and $v \in V_w$ there exists an unique morphism as before such that 
	$f(e) = v$.
\end{prop}
\textbf{Proof} 
\begin{enumerate}
	\item If $w=(a,b)$ and $a<b$. Let  $v = (v_1, v_2) \in  V_{(a,b)}$  and $f: P(0)\to V_0$ with $f(e_1) = v_1$;  $f(e_2) = v_2$.  Since $(0, v_1) \in  V_{(a,b)}$, then $f_w(0, e_1) = (0, v_1)$  and $f_w((e_1, e_2)) = (v_1, v_2) \in  V_{(a,b)}$; therefore  $f_w(P(w)_w) \in  V_w$.	
	Let $(a_1, b_1) \in \theta$, then $$P(w)_{(a_1,b_1)}=\langle(\lambda(a, a_1)e_1, 0), (0, \lambda(a, b_1)e_1), (\lambda(b, a_1)e_2, 0), (0, \lambda(b, b_1)e_2)\rangle.$$
	 If $\lambda(a, a_1)\neq 0,$
	then $a<a_1$ and therefore $(v_1, 0) \in  V_{(a_1,b_1)}$ and  $f_z(\lambda(a, a_1)e_1, 0)= (v_1, 0) \in  V_{(a,b)}$. In the same way,  it is seen that $f_z$, sends each generator
	from $P(w)_{(a_1, b_1)}$ into $V_{(a_1, b_1)}$. The uniqueness of $f$ is clear.
	
	\item If  $w = (a, b)$  and $a, b$ are incomparable. In this case, $P(w)_{(a,b)} =\langle
	(e_1, e_2) \rangle$  and $f_w((e_1, e_2)) = (v_1, v_2) \in  V_{(a,b)}$.  Therefore, $f_w(P(w)_{(a,b)}) \in 
	V_{(a,b)}$.  For the rest it is checked as in the previous case. 
	\item  If $w=\{a\}$, the proof  is similar to the previous cases.
\end{enumerate} 
\vspace{-0.4cm}
$\hfill\blacksquare$
\smallskip
\begin{prop}{\label{proyectives}} 
	The representations ${P}(w)=(P(w)_0,P(w)_z)_{z\in \theta}$ have the following properties:
	\begin{enumerate}
		\item ${P}(w)$ is an $\varepsilon$-projective representation.
		\item $\textnormal{End}({P}(w)) \cong k$ if  $w= (a, b)$ with $a$ and $b$ incomparable or when  $w$ consists of a single element.
		If $w=(a, b)$ with $a < b$ then $\textnormal{End}({P}(w)) \cong k[x]/x^2$.
		Therefore ${P}(w)$ is  indecomposable for all  $w \in \theta$.
		\item For any representation  $(V_0, V_z)_{z\in \theta}$, there exists an $\varepsilon$-deflation
		$g:(Q_0,Q_z)_{z\in \theta}\to (V_0,V_z)_{z\in \theta}$, where  $(Q_0,Q_z)_{z\in \theta}$ is $\varepsilon$-projective.
		\item If  $(Q_0,Q_z)_{z\in \theta}$ is an  indecomposable  projective representation  of $\rep(\mathscr{P},\theta)$, with $Q_z \neq 0$ for some $z\in \theta$, then $(Q_0,Q_z)_{z\in \theta}\cong {P}(w)$ for some $w\in \theta$.
	\end{enumerate}
\end{prop}
\textbf{Proof}\begin{enumerate} 	
	\item Let $f: (E_0, E_z)_{z\in \theta}\to (V_0, V_z)_{z\in \theta}$ an $\varepsilon$-deflation and $g: (P(w)_0, P(w)_z)_{z\in \theta}\to 	(V_0,V_z)_{z\in \theta}$ be a morphism. We take 
	$f^w(\underline{e}) \in  V_w$, where $\underline{e}$ is the generator of $(P(w)_0,P(w)_z)_{z\in \theta}$. Since $f^w$ is  surjective there exists $v_1\in E_w$ such that $f^w(v_1) = v$. By Proposition \ref{cubierta} there exists a morphism  $h: (P(w)_0,P(w)_z)_{z\in \theta}\to (E_0,E_z)_{z\in \theta}$ such that $h(\underline{e})=v_1$, so $fh(\underline{e}) = g(\underline{e})$. By the uniqueness in the Proposition \ref{cubierta}  is obtained that $fh=g$. Therefore,  $(P(w)_0,P(w)_z)_{z\in \theta}$  is a projective representation.
	
	\item  If $w=(a,b)$ or $w=\{a\}$ then $P(w)_w=\langle \underline{e}\rangle$ with $\underline{e}$ the generator of $(P(w)_0,P(w)_z)_{z\in \theta}$; therefore, if $f:(P(w)_0,P(w)_z)_{z\in \theta}\to (P(w)_0,P(w)_z)_{z\in \theta} $ then $f(\underline{e})=c\underline{e}$ with $c\in k$. Hence, $f=c(id_{{P(w)}_0})$.  This proves that 
	$$\textnormal{End}((P(w)_0, P(w)_z)_{z\in \theta} = k(id_{P(w)})
	\cong k.$$
	
	We suppose now that $w = (a, b)$ with $a < b$, then if $\underline{e} = (e_1, e_2)$ is the generator 
	of $P(w)_0$. We have that  $P(w)_w =\langle (e_1, e_2), (0, e_1) \rangle$.
	Let $f$ be an endomorphism of $(P(w)_0, P(w)_z)_{z\in \theta}$, then $f_w((e_1, e_2))=
	c(e_1, e_2)+d(0, e_1)$ with $c, d \in k$. Therefore $f(e_1)= ce_1, f(e_2) = ce_2+de_1$. In view of Proposition \ref{cubierta} the  morphism $f$ is completely determined by the matrix 
	$M(f) =\begin{pmatrix}
		c & d\\
		0 & c\\	
	\end{pmatrix}$. 
	If $f_1$ is another automorphism of $(P(w)_0,P(w)_z)_{z\in \theta}$, then $M(f_1f)=M(f_1)M(f)$. Hence,
	$$\textnormal{End}((P(w)_0,P(w)_z)_{z\in \theta})\cong 
	\left\lbrace\begin{pmatrix}
		c& d\\
		0 &c
	\end{pmatrix}\mid 
	c,d\in k \right\rbrace \cong k[x]/(x^2).$$
	
	\item For $V$ we choose a basis $B(0)$  and for each $z\in \theta$ such that $V_z\neq 0$, we choose $B(z)$ a $k$-basis of $V_z$. For each $v\in B(0)$ we take the morphism $f_v:{S}\to (V_0, V_z)_{z\in \theta}$
	which sends $1\in k$ in $v\in V$ and for $v \in  B(z)$ we have a morphism 
	$f_v:(P(w)_0, P(w)_z)_{z\in \theta}\to  (V_0, V_z)_{z\in \theta}$,  
	such that $f_v(\underline{e})=v$ where $\underline{e}$ is the  generator of $(P(w)_0, P(w)_z)_{z\in \theta}$. Let $B=\underset{z}\bigcup B(z)$, then we have a morphism 
	$$f = (f_v)_{v\in B}: \underset{v\in B(0)}\bigoplus {S}\quad \underset{z}\bigoplus \underset{v\in B(z)} \bigoplus {P}(w) \to (V_0, V_z)_{z\in \theta};$$
	clearly this morphism is an $\varepsilon$-deflation  and the representation  
	$$\underset{v\in B(0)}\bigoplus  {S}\underset{z}\bigoplus\underset{v\in B(z)}\bigoplus {P}(w),$$ is $\varepsilon$-projective.	
	
	\item Let $(Q_0,Q_z)_{z\in \theta}$ be a  projective representation,
	such that for some $z\in Q_z \neq 0$. From the above, we have a deflation:
	$${P}\xrightarrow{f}(Q_0,Q_z)_{z\in \theta},$$ 
	then there exists a morphism $h: (Q_0,Q_z)_{z\in \theta}\to {P}$ such that $fh = id_Q$. This implies that  $(Q_0,Q_z)_{z\in \theta}$ is a direct sum of  ${P}$. The last representation is a direct sum of representations ${S}$ and ${P}(z)$;  therefore, our representation is isomorphic
	to one of these, and as for some $z\in \theta$, $Q_z\neq 0$, then  $(Q_0,Q_z)_{z\in \theta} \cong {P}(w)$ for some $w\in \theta$.
\end{enumerate}
\vspace{-0.5cm}
$\hfill\blacksquare$
\smallskip
\begin{rem}
	An exact category is said to have enough projectives if it satisfies property 3 of Proposition $\ref{proyectives}$.
\end{rem}
\smallskip
\begin{defin}{\label{injec}}
	A representation $(I_0,I_z)_{ z\in \theta}$ is called $\varepsilon$-injective if given an $\varepsilon$-inflation
	$f:(V_0,V_z)_{z\in \theta} \to (E_0,E_z)_{z\in \theta} $ and a morphism $g:(V_0,V_z)_{z\in \theta}\to (I_0,I_z)_{z\in \theta}$, there exists a morphism $h:(E_0,E_z)_{z\in \theta}\to (I_0,I_z)_{z\in \theta}$  such that $hf=g$.
\end{defin}	
\smallskip
Henceforth it is convenient to use the following notation  to represent poset with an  involution:
the pair $(V_0,V_z)_{z\in \theta}$ where $V_0$ is a $k$-vector space and  $V_z\subset V_0^z$, is a representation of $\mathscr{(P,\theta)}$ if and only if for each  $x\leq y$	in $\mathscr{P}$, there exists a linear transformation $\tau: V_{[x]}\xrightarrow{} V_{[y]}$ such that  the following diagram is commutative 
\[  	 
\begin{tikzcd}
	V_{[x]}  \arrow[r, "i_{[x]}"]\arrow[d, "\tau"] &  V_0^{[x]} \arrow[d, "i_y\pi_x"]  \\
	V_{[y]}  \arrow[r, "i_{[y]}"] &  V_0^{[y]}  
\end{tikzcd}
\]

where $i_{[x]}:V_{[x]}  \to  V_0^{[x]}$ and $i_y:V_0\to V_0^{[y]}$ are the inclusions and $\pi_x:V_0^{[x]}\to V_0$ is the projection.

\subsection{Representations by quotients}

Let $(\mathscr{P},\theta)$ be a poset with an involution and $k$ be a field. 	A \textit{ representation  by quotient} $(V_0, j_z)_{z\in \theta}$, consists of a $k$-vector space  $V_0$ and for each $z\in \theta$ an epimorphism $j_z:V_0^{z}\to V_z$ such that if $a_1<a$ and $z=(a,b)$, $z_1=(a_1,b_1)$ then there exists a morphism $\tau: V_z\to V_{z_1}$ such that
$$\tau j_z=j_{z_1}i_{a_1}\pi_a.$$

A morphism $f:(V_0,j_z)_{z\in \theta}\to (V'_0,j'_z)_{z\in \theta}$
consists of a linear transformation  $f_0: V_0\to V'_0$  and for each $z\in \theta$ a linear transformation $f^z:V_z\to V'_z$ such that the following diagram commutes
\[  	 
\begin{tikzcd}
	V_0^z  \arrow[r, "j_z"]\arrow[d, "f_0^z"] &  V_z \arrow[d, "f_z"]  \\
	(V'_0)^z  \arrow[r, "j'_z"] &V'_z 
\end{tikzcd}
\]
We denote by  $\repq(\mathscr{P},\theta)$ the category of quotient representations.
\par\bigskip
\begin{prop}
	
	There are functors 
	$$C: \rep(\mathscr{P},\theta)\to \repq(\mathscr{P}^{op},\theta),\qquad \textnormal{and}\qquad K: \repq(\mathscr{P}^{op},\theta) \to \rep(\mathscr{P},\theta), $$
	
	defined by $C((V_0,V_z)_{z\in \theta})= (V_0, \textnormal{Coker}(i_z))_{z\in \theta}$, where $i_z:V_z\to V_0^{z}$ is the inclusion and 
	$K((V_0, j_z))_{z\in \theta}= (V_0, \textnormal{Ker}(j_z))_{z\in \theta}$. 
	Further,  $CK\cong id_{\repq(\mathscr{P}^{op},\theta)}$ and $KC\cong id_{\rep(\mathscr{P},\theta)}$
	therefore  $\rep(\mathscr{P},\theta)$ is equivalent to  $\repq(\mathscr{P}^{op},\theta)$.
\end{prop} 
\textbf{Proof} Let $(V_0, V_z)_{z\in \theta}$  be an object
of $\rep(\mathscr{P},\theta)$
and we take $j_z:V_0^z\to V_z'$ the cokernel of $i_z$. We suppose that $x\in z$ and $y\in z_1$ with $x<y$ then we obtain the morphism $i_y\pi_x:V_0^{z}\to V_0^{z_1}$ and a morphism $\tau: V_z\to V_{z_1}$ such that $i_{z_1}\tau=i_y\pi_x i_z$. Therefore there exists a morphism $\tau':V'_z\to V'_{z_1}$ such that  the following diagram is commutative
\[ A:
\begin{tikzcd}
	V_z  \arrow[r, "i_z"]\arrow[d, "\tau"] &  V_0^{z} \arrow[d, "i_y\pi_x"] \arrow[r, "j_z"]  
	&V'_z \arrow[d, "\tau'"] \\
	V_{z_1}  \arrow[r, "i_{z_1}"] &  V_0^{z_1} \arrow[r, "j_{z_1}"]  &
	V_{z_1}' 
\end{tikzcd}\]

This proves that $(V_0, j_z)_{z\in \theta}\in \repq(\mathscr{P}^{op},\theta)$.

\par\bigskip

Now, let $f:(V_0,V_z)_{z\in \theta}\to (W_0,W_z)_{z\in \theta}$ be a morphism in $\rep(\mathscr{P},\theta)$; we denote  $r_z:W_z\to W_0^{z}$ the inclusion  and by $r_z':W_0^{z}\to W_z'$ its cokernel. The morphism $g_z:V_z\to W_z$ is obtained, and it is  such that $f_0^{z}i_z=r_zg_z$. Therefore there exists a morphism $f_z:V_z'\to W_z'$  such that the following diagram is commutative
\[ B:  	 
\begin{tikzcd}
	V_z  \arrow[r, "i_z"]\arrow[d, "g_z"] &  V_0^{z} \arrow[d, "f_0^z"] \arrow[r, "i_z'"]  
	&V'_z \arrow[d, "f_z"] \\
	W_{z}  \arrow[r, "r_{z}"] &  W_0^{z} \arrow[r, "r_{z}' "]  &
	W_{z}' 
\end{tikzcd}
\]
thus,  $f_0^{z}$ is a morphism of $C((V_0,V_z)_{z\in \theta})$ in $C((W_0,W_z)_{z\in \theta})$. We define $C(f)=f_0^{z}$.

\par\bigskip
Now, if $(V_0,j_z)_{z\in \theta}$ is an object of $\repq(\mathscr{P}^{op},\theta)$, by using  diagram A,  is obtained that  $K((V_0,j_z)_{z\in\theta})\in  \rep(\mathscr{P},\theta)$. If $f:(V_0,j_z)_{z\in\theta}\to (W_0,r_z')_{z\in\theta}$ is a morphism in $\repq(\mathscr{P}^{op},\theta)$ such that $f_0:V_0\to W_0$ then by using B is obtained that $f_0$ produces a morphism of $K((V_0,j_z)_{z\in\theta})$ in $K((W_0,r_z')_{z\in\theta})$. The rest of the proof  is clear. 
$\hfill\blacksquare$

\bigskip
Henceforth,  if $W$ is a $k$-vector space $D(W)=\text{Hom}_k(W,k).$ 
\smallskip
\begin{prop}
	There are contravariant functors
	$$D_1:\rep(\mathscr{P},\theta)\to \repq(\mathscr{P},\theta), \qquad \textnormal{and} \qquad D_2:\repq(\mathscr{P},\theta) \to \rep(\mathscr{P},\theta),$$
	defined by  $D_1((V_0, V_z)_{z\in \theta})= (D(V_0),D(i_z))_{z\in \theta}$ where $i_z:V_z\to V^z$ is the inclusion  and $D(i_z):D(V^z)=D(V)^z\to D(V_z)$ 
	and $D_2((V, j_z)_{z\in \theta})= (D(V),\text{im}(D(j_z)))$. Further  $D_2D_1\cong id_{\rep (\mathscr{P},\theta)}$ and $D_1D_2\cong \text{id}_{\repq (\mathscr{P}, \theta)}$.
\end{prop}
\textbf{Proof.} We identify $D(V_0^{z})=D(V_0)^z$. Let 
$(V_0, V_z)_{z\in \theta} \in  \rep(\mathscr{P},\theta)$, then $D_1((V_0, V_z)_{z\in \theta}) = (D(V_0),D(i_z))_{z\in \theta}$ where $i_z: V_z \to V_0^{z}$ is the inclusion. Then if  $a \in z$, $a_1 \in z_1$ with $a_1 < a$. Hence there exists a morphism  $\tau: V_{z_1} \to  V_z$ such that $i_z\tau = i_a\pi_{a_1}i_{z_1};$
therefore, 
$$D(\tau)D(i_z) = D(i_{z_1})D(\pi_{a_1}) D(i_a),$$
We observe that  $D(i_{a_1}): D(V_0)^{z_1} \to  D(V_0)$ is equal to $\pi_{a_1}$ and $D(\pi_a): D(V_0) \to D(V_0)^z$ is equal to $i_a$; therefore
$$D(\tau)D(i_z)= D(i_{z_1})i_{a_1}\pi_{a}.$$

The above implies that $D_1((V_0, V_z)_{z\in \theta}) \in  \rep(\mathscr{P},\theta)$. It is clear that 

$$f: (V_{0},V_{z})_{ z \in \theta} \to (W_{0},W_{z})_{z\in \theta}$$

is a morphism in $\rep(\mathscr{P},\theta)$, then $D(f_0):D(W_0) \to D(V_0)$ determines a morphism $D_1(f): D_1((W_{0},W_{z})_{z\in \theta}) \to D_1((V_0, V_z)_{z\in \theta})$.
The rest of the proposition proceeds in a similar way.
$\hfill\blacksquare$
\smallskip
\begin{defin}
	We consider $\varepsilon_q$ the class of sequences in $\repq (\mathscr{P}, \theta)$
	which have  the form 
	\[
	(V^{1}_0,j^{1}_z)_{z\in\theta}\xrightarrow{f} (V^{2}_0,j^{2}_z)_{z\in\theta} \xrightarrow{g} (V^{3}_0,j^{3}_z)_{z\in\theta},	
	\]
	such that 
	\[
	0 \xrightarrow{} V^{1}_0 \xrightarrow{f_0} V^{2}_0 \xrightarrow{g_0} V^{3}_0	\xrightarrow{} 0,
	\] and 
	\[
	0 \xrightarrow{} V^{1}_z \xrightarrow{f_z} V^{2}_z \xrightarrow{g_z} V^{3}_z	\xrightarrow{} 0,
	\]
	are exact, where $j^{i}_z: (V^{i}_0)^z\to V^{i}_z$.
\end{defin}
\smallskip
\begin{prop}
	The functor $D_1$ sends $\varepsilon$-sequences  to $\varepsilon_q$-sequences and the functor $D_2$ sends $\varepsilon_q$-sequences in $\varepsilon$-sequences. In particular, a morphism $f:(V_0, V_z)_{z\in \theta}\to (V'_0, V'_z)_{z\in \theta}$ in  $\rep (\mathscr{P}, \theta)$ is an $\varepsilon$-inflation ($\varepsilon$-deflation, respectively) if and only if $D_1(f)$ is an $\varepsilon$-deflation ($\varepsilon$-inflation, respectively).
\end{prop}	
\smallskip
\begin{corol}
	The class of morphisms $\varepsilon_q$ is an exact structure. Further the category $\repq (\mathscr{P}, \theta)$ has enough injectives. 
\end{corol}
\smallskip
\begin{corol}
	The exact category  $(\rep (\mathscr{P}, \theta), \varepsilon)$ has enough injectives. The indecomposable injectives of this category have the form   $KD_1(P_z)$ for $z\in \theta$ and $KD_1((k,0_z)_{z\in \theta})$, where $P_z$ and $(k, 0_z)_{z\in \theta}$ are projectives in $\rep(\mathscr{P}^{op}, \theta).$
\end{corol}
\textbf{Proof} The indecomposable injectives of $\repq (\mathscr{P}^{op}, \theta)$ have  the form $D_1(P_z)$ and $D_1((k,0_z))_{z\in \theta}$. Since the functor $K$ is an equivalence of categories such that sends $\varepsilon_q$-sequences in $\varepsilon$-sequences then the injectives indecomposables of $\rep (\mathscr{P}, \theta)$ are the form $KD_1(P_z)$ for $z\in  \theta$ and $KD_1((k, 0_z)_{z\in \theta})$.
$\hfill\blacksquare$
\smallskip
\begin{rem}
	$KD_1((k, 0_z)_{z\in \theta})=J$ is the representation  $J=(J_0,J_z)_{z\in \theta}$ such that $J_0=k$ and $J_z=k^z$.
\end{rem}

\section{The endomorphism algebra}{\label{endomorphism}}
Let  $(\mathscr{P}, \theta)$  be a poset with an involution. We know that in the exact category  $(\rep(\mathscr{P}, \theta), \varepsilon)$ a system of representatives of  isomorphism classes 
of indecomposable projectives
is given by ${P}(z)$ for $z\in \theta$ and ${P}(0)=S$, as in Remark 2. 
We take  
$\underline{P}=\underset{z\in \theta}\bigoplus{P}(z)\bigoplus P(0)$  and   $A=\text{End}_{\rep(\mathscr{P}, \theta)}(\underline{P})$. Take    $P(j)$ with $j = z \in \theta$ or
$j = 0$. If we consider the projection $\pi_j: \underline{P} \to P(j)$ and the inclusion $\sigma_j : P(j) \to \underline{P}$ we
obtain the idempotent $e_j = \sigma\pi_j \in  A$. Then $1_A=\underset{z\in \theta}\sum e_z+e_0,\qquad \textnormal{and}\qquad
A=\underset{z\in \theta}\bigoplus e_zA\bigoplus e_0A.$\\
Now we recall that an algebra $B$ is called right-peak algebra if $\textnormal{soc}(B_B)$ is a sum of
copies of a simple projective module.

\begin{prop}
	The algebra $A = \textnormal{End}_{\rep(\mathscr{P}, \theta)}(\underline{P})$ is a right peak algebra and the only
	simple projective right A-module up to isomorphism is $e_0A$.
\end{prop}
\textbf{Proof}
We have:
$$A =\begin{pmatrix}
	A_1& M\\
	0 & k
\end{pmatrix}$$

where $A_1=\textnormal{End}_{\rep(\mathscr{P},\theta)}(\underset{z\in \theta}\bigoplus P(z))$ and $M = \textnormal{Hom}_{\rep(\mathscr{P}, \theta)}(P(0), \underset{z\in \theta}\bigoplus P(z))$. Here $M$ is
a faithful left $A$-module so by \cite{Simson}, $A$ is a right peak algebra and the only simple
projective $A$-module up to isomorphism is $e_0A$.
$\hfill\blacksquare$
\begin{rem}{\label{rema}}
	The injective envelope of $e_0A$ is $E = D(Ae_0)$ and
	$$\textnormal{dim}_kEe_z = \textnormal{card}(z),\quad  \textnormal{dim}_kEe_0 = 1$$
\end{rem}
\textbf{Proof} Here
{\footnotesize{
		\begin{equation*}
			\begin{split}
				\textnormal{dim}_kEe_z = \textnormal{dim}_kD(Ae_0)e_z
				= \textnormal{dim}_kD(e_zAe_0)
				= \textnormal{dim}_ke_zAe_0
				= \textnormal{dim}_k\textnormal{Hom}_{\rep(\mathscr{P}, \theta)}(P(0), P(z))
				= \textnormal{card}(z).
			\end{split}
\end{equation*}}}
Now, we consider the full and faithful functor $H = \textnormal{Hom}_{\rep(\mathscr{P}, \theta)}(\underline{P}, -): \rep(\mathscr{P}, \theta)\to 
\textnormal{mod}A$. Since $\underline{P}$ is projective, this functor sends $\varepsilon$-sequences in exact sequences in $\textnormal{mod}A$.
Moreover, $H(\underline{P}) = A$ and $H(P(z)) = e_zA$.
$\hfill\blacksquare$

\begin{rem} 
	If $f: V \to W$ is a morphism in $\rep(\mathscr{P}, \theta)$, with $f:V_0 \to W_0$ a monomorphism
	then $H(f): H(V ) \to H(W)$ is also a monomorphism.
\end{rem}

\begin{lema}
	If $V \in \rep(\mathscr{P}, \theta)$, then $\textnormal{soc}H(V )$ is projective. Morerove,  $\textnormal{dim}_k(\textnormal{soc}H(V )) =
	\textnormal{dim}_kV_0$.
\end{lema}
\textbf{Proof} Take $V \in \rep(\mathscr{P}, \theta)$ and consider $\textnormal{soc}H(V)$, if this is not a direct sum of
copies of $e_0A$, there is a non zero morphism $u: e_zA \to \textnormal{soc}(H(V ))$ with $z \in \theta$, such
that $u(\textnormal{rad}e_zA) = 0$. We have $u = H(u')$ with $u': P(z) \to V$, so $u': P(z)_0 \to V_0$ is a non-zero $k$-linear map, then there is a non zero $k$-linear map $\lambda: {S}_0 \to  P(z)_0$
such that $u'\lambda\neq 0$, but $\lambda$ is a morphism ${S} \to P(z)$ then $H(u')H(\lambda) \neq 0$. Here the
image of $H(u')$ lies in the radical of $e_zA$, so $H(u')H(\lambda) = u(H(\lambda)) = 0$ a contradiction. This implies that $\textnormal{soc}H(V )$ is projective. Then $H(V ) \in  \textnormal{mod}_{\textnormal{sp}}(A)$. Moreover,
$\textnormal{dim}_k(\textnormal{soc}H(V)) = \textnormal{dim}_k\textnormal{Hom}_A(e_0A, H(V)) = \textnormal{dim}_k\textnormal{Hom}_{\rep(\mathscr{P}, \theta)}({S}, V) = \textnormal{dim}_kV_0$.
$\hfill\blacksquare$

\begin{lema} 
	Consider $J \in  \rep(\mathscr{P}, \theta)$, the representation of \textnormal{Remark 5} , then $H(J) \cong E$ the
	injective envelope of $e_0A$.
\end{lema}
\textbf{Proof} Since $J_0 = k$, then $\textnormal{soc}H(J)$ is a simple right $A$-module, the injective envelope
of $H(J)$ is of the form $u: H(J) \to E$. Here $\textnormal{dim}_kH(J)e_z = \textnormal{dim}_k\textnormal{Hom}_{\rep(\mathscr{P}, \theta)}(P(z), J)=
\textnormal{dim}_k(J_z) = \textnormal{dim}_k(k^z) = \textnormal{card}(z)$. But because the above $\textnormal{dim}_kEe_z = \textnormal{dim}_kH(J)e_z, $ therefore $H(J)\cong E$.
$\hfill\blacksquare$

\begin{teor}{\label{teorem}}
	The functor $H$ induces an equivalence of categories $H : \rep(\mathscr{P}, \theta)\to \textnormal{mod}_{\textnormal{sp}}(A)$ where $\textnormal{mod}_{\textnormal{sp}}(A)$ is the full subcategory of $\textnormal{mod}A$ whose objects are the right
	$A$-modules $M$ with $\textnormal{soc}M$ projective.
\end{teor}
\textbf{Proof} We know that $H$ induces a full and faithful functor $H: \rep(\mathscr{P}, \theta)\to  \textnormal{mod}_{\textnormal{sp}(A)}$.
Now we are going to prove that the functor $H$ is dense. For this take $X \in \textnormal{mod}_{\textnormal{sp}}A$, then
the injective envelope of $X$ is of the form $X \to E^n \cong H(J^n)$ so we have a monomorphism
$f : X \to  H(J^n)$. Take $g : Q\to  X$ a projective cover of $X$ then we may assume $Q = H(L)$
with $L$ a projective in $\rep(\mathscr{P}, \theta)$. In this way we obtain a morphism $gf: H(L) \to H(J^n
)$, since $H$ is a full functor there is a morphism $h_1: L \to J_n$ with $H(h)= gf$. The morphism
$h = h_2h_1$ with $h_1 : L\to  h(L)$ a deflation and $h_2: h(L)\to  J^n$
such that $h_2 : h(L)_0 \to  (J^n)$
a monomorphism. Therefore $gf = H(h) = H(h_2)H(h_1)$ where $H(h_1)$ an epimorphism and
$H(h_2)$ a monomorphism. From this we infer that $\textnormal{Im}(gf) = \textnormal{Im}(f) \cong \textnormal{Im}(H(h_1))$. Here
$\textnormal{Im}(f) = X$ and $\textnormal{Im}(h_1) = H(\textnormal{Im}h_1)$. This proves our proposition.
$\hfill\blacksquare$

\begin{rem}
	By a result of D. Simson \textnormal{[see, \cite{Simson}]}, the category $\textnormal{mod}_{\textnormal{sp}}(A)$ has almost split sequences. Then \textnormal{Theorem} $\ref{teorem}$, implies that the exact category $(\rep(\mathscr{P}, \theta),\varepsilon )$ has also almost
	split sequences.
\end{rem}
\section{The
	Auslander-Reiten quiver of posets with  an involution of type $\mathfrak{U}_n$.}

By using the results from the previous section, we construct the Auslander-Reiten quiver for a poset type that we will denote by  $ \mathfrak{U}_n $.  In this section we will assume that $k$ is an algebraically closed field. 

\subsection{Poset with an involution of Type $\mathfrak{U}_n$}{\label{sect5.1}}
We denote by $\mathfrak{U}_n$ to the poset with an  involution $(\mathscr{P}, \leq, \theta)$ where  
$(\mathscr{P}, \leq)=\{a_n<a_{n-1}<\cdots<a_{1}<b_1<b_2<\cdots <b_{n-1}<b_{n}\}$ and $\theta= \{(a_i,b_i)\}_{i=1,\dots,n}$.  We denote by $\textnormal{Rep}(\mathfrak{U}_n)$ the category of representations of the poset $\mathfrak{U}_n$. 
\par\bigskip

The Hasse diagram of the poset $\mathfrak{U}_n$ is as follows
\begin{figure}[h]
	\setlength{\unitlength}{1pt}
	\begin{center}
		
		\begin{center}
			\begin{tikzpicture}
				
				\node (1) at (-2,-1.5) {$(\mathscr{P},\theta)$=};
				\node (1) at (-0.3,0.11) {$b_n$};
				\node (1) at (-0.5,-0.6) {$b_{n-1}$};
				\node (1) at (-0.3,-1.2) {$b_1$};
				\node (1) at (-0.3,-1.8) {$a_{1}$};
				\node (1) at (-0.5 ,-2.4) {$a_{n-1}$};
				\node (1) at (-0.4,-3) {$a_{n}$};
				
				\node (0) at (0,0) {};
				\node (1) at (0,0) {$\bullet$};
				\node [below  of=1, node distance=0.6cm] (1')  {$\bullet$};
				\node [below  of=1', node distance=0.6cm] (1'')  {$\bullet$};
				\node [below  of=1'',node distance=0.6cm] (2')  {$\bullet$};
				\node (1) at (0,-2) {$\vdots$};
				\node [below  of=2',node distance=0.6cm] (3')  {$\bullet$};
				\node (1) at (0,-0.8) {$\vdots$};
				\node [below  of=3',node distance=0.6cm] (4')  {$\bullet$};
				
				\draw [  shorten <=-2pt, shorten >=-2pt] (0) -- (1);
				\draw [  shorten <=-2pt, shorten >=-2pt] (1'') -- (2');
				\draw [ shorten <=-2pt, shorten >=-2pt] (3') -- (4');
			\end{tikzpicture}
		\end{center}
	\end{center}
	\caption{Hasse diagram  of poset with an involution of type $\mathfrak{U}_n$ }
\end{figure}
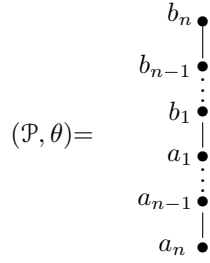

\smallskip
We consider the following representations of $\mathfrak{U}_n:$

\begin{description}
	\item[ a.] $\mathcal{L}_{1,i}= (\mathcal{L}_0, \mathcal{L}_{(a_j,b_j)})_{j\geq 1}$
	where $\mathcal{L}_0=k\{e\}$ and 
	\begin{align}{\label{l1i}}
		\mathcal{L}_{(a_j,b_j)}&=\begin{cases}
			(0,0), &\quad \quad \text{if } \quad  j<i,\\
			\langle(0,e)\rangle,&\quad  \quad  \text{if } \quad j\geq i.\\
		\end{cases}
	\end{align}
	
	\item[b.] $\mathcal{L}_{2,i}= (\mathcal{L}_0, \mathcal{L}_{(a_j,b_j)})_{j\geq 1}$
	where $\mathcal{L}_0=k\{e\}$ and 
	\begin{align}\label{l2i}
		\mathcal{L}_{(a_j,b_j)}&=\begin{cases}
			\langle(0,e),(e,0)\rangle, & \quad \text{if } \quad  j\leq i, \\
			\langle(0,e)\rangle, & \quad \text{if } \quad j> i.\\
		\end{cases}
	\end{align}

	\item[c.] $\mathcal{L}_{3,i}= (\mathcal{L}_0, \mathcal{L}_{(a_j,b_j)})_{j\geq 1}$
	where 	$\mathcal{L}_0=k\{e_1,e_2\}$ and 
	\begin{align}{\label{l3i}}
		\mathcal{L}_{(a_j,b_j)})&=\begin{cases}
			\langle (0,e_1),(e_1,0)\rangle,& \quad \text{if } \quad j<i, \\
			\langle (0,e_1),(e_1,e_2)\rangle,& \quad\text{if } \quad  j=i,\\
			\langle (0,e_1),(0,e_2)\rangle,&  \quad \text{if } \quad j> i.\\
		\end{cases}
	\end{align}
\end{description}

It is clear that each one of the previous representations are indecomposable.
\smallskip
\begin{rem}
	The representation ${S}=(k, S_z)_{z\in \theta}$
	with $S_z=0$ for all $z\in \theta$, is called a trivial indecomposable  representation of  $\mathfrak{U}_n$. 
\end{rem}
\smallskip
\begin{prop}{\label{cadena1}}
	The representations above is the complete list of non trivial indecomposable representations  of $\mathfrak{U}_n$.
\end{prop}
\textbf{Proof.} 
We will prove by induction on $n$ that any representation of $\mathfrak{U}_n$ can be written as a
direct sum of some representations which are isomorphic to the trivial representation and
to some representations in the previous list.
We first consider the case $n = 1$, in this case $\mathfrak{U}_1 = (\mathscr{P}, \theta)$, with $P = \{a, b\}$ with $a < b$
and $\theta= z$, where $z = (a, b)$. Take $V = (V_0, V_z)$ a representation of $\mathfrak{U}_1$. Here we can take $V_0^{z}= V_0 \oplus V_0$ where for $w = (w_1, w_2) \in V_0^{z}$, $\pi_a(w) = w_1$, and $ \pi_b(w) = w_2$. Then by definition if
$(w_1, w_2) \in  V_z,(0, w_1) \in  V_z$. Take $u(1), \dots , u(n)$ a basis for $V_z$ and $e_1,\dots, e_m$ a basis for $V_0$,
then we have $u(s) = (\underset{j=1}{\overset{m}\sum} \alpha_{j,s}^1 e_j, \underset{j=1}{\overset{m}\sum} \alpha_{j,s}^2 e_j)$ for $ s = 1,\dots, n$. We have the matrices
$T_a = (\alpha_{j,1}^1)$ and $T_b = (\alpha^2_{j,i})$
Choosing some other bases, the matrices $T_a$ and $T_b$ change to $ST_aR$ and $ST_bR$ where
$S$ and $R$ are non singular matrices.
Here we are assuming that $k$ is an algebraically closed field, then using the Kronecker
decomposition theorem we may find bases for $V_z$ and $V_0$ such that the corresponding
matrices $T_a, T_b$ with respect to the chosen bases is a direct sum of matrices of the form:

\begin{multicols}{3}
	\begin{itemize}
		\item[\textbf{a.}]  $\begin{pmatrix}
			E_n\\
			0_{1,n}
		\end{pmatrix},\begin{pmatrix}
			0_{1,n}	
			\\E_n
		\end{pmatrix}$
		
		\item[\textbf{b.}]  $\begin{pmatrix}
			0_{1,n}\\
			E_n
		\end{pmatrix},\begin{pmatrix}
			E_n\\
			0_{1,n}	
		\end{pmatrix}$
		\item[\textbf{c.}] $(0_{n,1},E_n), (E_n,0_{n,1})$\\
		\item[\textbf{d.}] $(E_n,0_{n,1}), (0_{n,1},E_n)$ 
		\item[\textbf{e.}] $J_{\lambda,n},E_n$
		\item[\textbf{f.}] $E_n, J_{\lambda,n}$ 
	\end{itemize}
\end{multicols}

Therefore we have a decomposition $V_0 = V^{(1)}_0\oplus \dots \oplus V^{(l)}_0$
and $V_z = V^{(1)}_z \oplus \dots \oplus  V^{(l)}_z$
with $V^{(j)}_z \subset (V^{(j)}_0)^z$
for $j = 1,\dots, l$. Moreover,  there are basis for each $V^{(j)}_z$ and $V^{(j)}_0$
such
that the corresponding matrices $T^j_{a} , T^j_{b}$ have one of the forms \textbf{a}, \textbf{b}, \textbf{c}, \textbf{d}, \textbf{e}, \textbf{f}. Since $(V_0, V_z)$ is
a representation, each $(V^{(j)}_{0}, V^{(j)}_{z} )$ is a representation of $(\mathscr{P}, \theta)$.\\

Therefore we may assume that there is a basis for $V_0$ and $V_z$ such that the pair of
matrices $T_a, T_b$ have one of the forms \textbf{a}, \textbf{b}, \textbf{c}, \textbf{d}, \textbf{e}, \textbf{f}. Here $\pi_a(V_z) \subset  \pi_b(V_z)$, then the only pair of matrices satisfying this condition are  \textbf{c}, \textbf{d}, \textbf{e}, \textbf{f}.\\

Consider the case \textbf{c}:\\
Here for $n = 1, V_z$ is generated by the vectors $u(1) = (0, e_1), u(2) = (e_1, 0)$, therefore
$(V_0, V_z) = \mathcal{L}_{2,1}$. Now if $n \geq 2$, then $V_z$ is generated by the vectors $u(1) = (0, e_1), u(2) =
(e_1, e_2), u(3) = (e_2, e_3), ..., u(n) = (e_{n-1}, e_n), u(n + 1) = (e_n, 0)$.
Observe that for $n = 2, (e_2, 0) \in V_z$ but the vector $(0, e_2)$ is not in $V_z$, for $n \geq 3$, the
vector $(e_2, e_3) \in V_z$ but $(0, e_2)$ is not in $V_z$. Therefore the only posibility is $n = 1$ and in this case $(V_0, V_z) =\mathcal{L}_{2,1}$. Similarly in case \textbf{d}, $n = 1$ and $(V_0, V_z) = \mathcal{L}_{2,1}.$\\

In cases \textbf{e} and \textbf{f} if $\lambda \neq 0$, the vectors in $V_z$ have the form $(\lambda x, x)$ in case \textbf{e} and $(x, \lambda x)$
in case \textbf{f}. In the first case for $x \neq 0,\lambda x \in \pi_a(V_z)$ but $(0, \lambda x)$ is not in $V_z$; similarly in the
second case there is $x \in \pi_a(V_z)$ but $(0, x)$ is not in $V_z$.
Suppose now $\lambda = 0$. In case \textbf{e}, for $n = 1, V_z =\langle (0, e_1)\rangle$, therefore $V = (V_0, V_z) = \mathcal{L}_{1,1}$.
For $n = 2, V_z =\langle(0, e_1),(e_1, e_2) \rangle$, in this case $V = (V_0, V_z) = \mathcal{L}_{3,1}$. For $n > 2$, we have
$(0, e_1),(e_1, e_2),(e_2, e_3) \in V_z$, but $(0, e_1)$ is not in $V_z$. Therefore this case can not happen.
For \textbf{f} with $\lambda = 0$, we have that $(e_1, 0)$ is in $V_z$ but $(0, e_1$) is not in $V_z$, therefore this
case can not happen. We have proved that the indecomposable representations of $(\mathscr{P},\theta)$ are
$\mathcal{L}_{1,1}, \mathcal{L}_{2,1}$ and $\mathcal{L}_{3,1}$.
This shows our result for case $n=1$.\\

Now, assume our result true for $u_{n-1}$. Let $V= (V_0, V_{(a_j ,b_j)})_{j\leq n}$ be a representation of
$\mathfrak{U}_n$. We take $\mathfrak{U}_1 = (a_n< b_n; \{a_n, b_n\})$. Then $(V_0, V_{(a_n,b_n)})$ is a representation of $\mathfrak{U}_1$.
Therefore $(V_0, V_{(a_n,b_n)}) = S(W_0)\oplus \tilde{V}$  with $\tilde{V}_0= V_1\oplus V_2\oplus V_3\oplus V_3'$
and
$$\tilde{V}_{(a_n,b_n)} =V_{(a_n,b_n)} = (0, V_1)\oplus (V_2, 0)\oplus (0, V_2)\oplus (0, V_3)\oplus H_\phi(V_3, V_3')$$ where, $V_0 = W_0\oplus V_1\oplus V_2\oplus V_3\oplus V_3', \phi: V_3\to  V_3'$ is an isomorphism and $H_\phi(V_3, V_3') = \{(w, \phi(w))\mid w\in  V_3\}$.\\

Let $\mathfrak{U}_{n-1}=\{a_{n-1}<\dots<a_1<b_1<\dots <b_{n-1}; (a_1,b_1), \dots,(a_{n-1},b_{n-1})\}$
and $\underline{V}=(V_0,V_{(a_i,b_i)})_{i\leq n-1}$ be the restriction of the representation $V$ to $\mathfrak{U}_{n-1}$.
\smallskip
\begin{rem}{\label{rem1}}
	If $u\in V_2\bigoplus V_3$ then $(u,0)$ and $(0,u)$
	are in $V_{(a_i,b_i)}$ for all $i<n$. Indeed, if $u\in V_2$, $(u,0)\in V_{(a_n,b_n)}$. As $a_n<a_i$, $(u,0)\in V_{(a_i,b_i)}$ and   as $a_i<b_i$, $(0,u)\in V_{(a_i,b_i)}$ too. If $u\in V_3$ then $(u,\phi(u))\in V_{(a_n,b_n)}$. As $a_n<a_i$, $(u,0)\in V_{(a_i,b_i)}$;  therefore 
	$(0,u)\in V_{(a_i,b_i)}$.
\end{rem}

\begin{rem}{\label{rem2}}
	If $(u,v)\in V_{(a_i,b_i)}$ with $i<n$ then $u$ and $v$ are in $V_1\bigoplus V_2\bigoplus V_3$. Indeed, as $a_i<b_n$, $(0,u)\in (a_n,b_n)$, therefore $u\in V_1\bigoplus V_2 \bigoplus V_3$.
	Analogously, as $b_i<b_n$ then $(0,v)\in V_{(a_n,b_n)}$, thus $v\in  V_1\bigoplus V_2 \bigoplus V_3$.
\end{rem}
\smallskip
We consider $V'= (V_1, V'_{(a_i,b_i)})_{1\leq i\leq n-1}$, where $V'_{(a_i,b_i)}= \{(u, v) \in  V_{(a_i,b_i)}\mid u, v \in  V_1\}$, and
$$V''= (L, V''_{(a_i,b_i)})_{1\leq i\leq n-1},$$ with $L= V_2\oplus V_3\oplus V'_3\oplus W'$ and
$V''_{(a_i,b_i)} = \{(u, v)\in  V_{(a_i,b_i)}\mid u, v\in L\}$. By Remarks \ref{rem1} and \ref{rem2}, for all $i$,   $1 \leq i \leq n-1$:
$$V''_{(a_i,b_i)} = (V_2, 0) + (0, V_2) + (V_3, 0) + (0, V_3).$$

Both $V'$ and $V''$ are representations of $\mathfrak{U}_{n-1}$.\\

\textbf{Claim} $V = V'\oplus  V''$.

We just have to prove that
for all $i<n$,
$$V_{(a_i,b_i)}=V'_{(a_i,b_i)}\bigoplus V''_{(a_i,b_i)}.$$

Let $(u,v)\in V_{(a_i,b_i)}$, due to Remark $\ref{rem2}$, $u=u_1+u_2$, $v=v_1+v_2$ where $v_1,v_2$ are in $V_1$, and $u_2,v_2$ are in $V_2\bigoplus V_3$; thus $(u_2,v_2)=(u_2,0)+(0,v_2)\in V_{(a_i,b_i)}$, so $(u_1,v_1)\in V_{(a_i,b_i)}$ and since $u_1,v_1$ are in $V_1$,  it is obtained that $(u_1,v_1)\in V'_{(a_1,b_1)}$ and  clearly $(u_2,v_2)\in V''_{(a_i,b_i)}$. This proves our claim for $n = 1$. We now suppose proved our result for $\mathfrak{U}_{n-1}$, for proving it for $\mathfrak{U}_n$ we denote by $\mathcal{L}^{n-1}_{1,s}, \mathcal{L}^{n-1}_{2,s}$ and $\mathcal{L}^{n-1}_{3,s}$ the representations given in (1), (2)
and (3) for $1 \leq  s \leq n-1$. In a similar way, we have the corresponding representations for
$u_n$. $\mathcal{L}^{n}_{1,t},\mathcal{L}^n_{2,t},\mathcal{L}_{3,t}$ for $1 \leq t \leq n$.
Then observe that for $1 \leq s, i\leq n-1$ one has:
\vspace*{-0.5cm}
\begin{multicols}{2}
	\begin{center}
		\begin{align*}
			\qquad \qquad \mathcal{L}^{n-1}_{1,s}(e)(a_i,b_i)&= \mathcal{L}^n_{1,s}(e)(a_i,b_i); &\qquad \mathcal{L}^n_{1,s}(e)(a_n,b_n) &=<(0, e)>;\\ \mathcal{L}^{n-1}_{2,s}(e)(a_i,b_i)&= \mathcal{L}^n_{2,s}(e)(a_i,b_i);&\qquad \mathcal{L}^n_{2,s}(e)(a_n,b_n)&=<(0, e)>;\\
			\mathcal{L}^{n-1}_{3,s}(e,f)(a_i,b_i)&= \mathcal{L}^n_{3,s}(e, f)(a_i,b_i); &\qquad \mathcal{L}^n_{3,s}(e, f)(a_n,b_n)&=<(0, e),(0, f)>.
		\end{align*}
	\end{center}
\end{multicols}

By application of the induction hypothesis we obtain:
{\footnotesize{
$$V'=\left(\underset{1\leq s\leq n-1}{\underset{e\in B_{1,s}}\bigoplus} \mathscr{L}^{1,s}_{n-1}(e)\right) \bigoplus  \left(\underset{1\leq t\leq n-1}{\underset{e\in B_{2,t}}\bigoplus} \mathscr{L}^{2,t}_{n-1} (e) \right)\bigoplus \left(\underset{1\leq r\leq n-1}{\underset{e\in B_{3,r}}\bigoplus} \mathscr{L}^{1,r}_{n-1} (e,\psi(e))\right)\bigoplus S(W'_0)$$}}

where $B_{1,s}, 1 \leq s \leq n-1, B_{2,t}, 1 \leq t \leq n-1, B_{3,r}, 1 \leq r \leq n-1, B'_{3,r}, 1 \leq r \leq n-1$ are
subsets of $V_1, \psi: B_{3,r}\to B'_{3,r}$ is a bijection, and $B= \underset{s}\bigcup B_{1,s}\, \underset{t} \bigcup  B_{2,t}\, \underset{t}\bigcup B_{3,t}\underset{r}\bigcup B'_{3,r}$ is a set of linearly independent elements and $V_1 =\langle B \rangle \oplus W'_0$.\\

Take now $h_1,\dots, h_l$ a $k$-basis for $V_2$ and $g_1,\dots, g_m$ a $k$-bases for $V_3$. Then for any $i$ with
$1 \leq  i \leq n-1$, we have
{\footnotesize{
\begin{align*}
	V(a_i,b_i)&=V'(a_i,b_i)\oplus V''(a_i,b_i)\\
	&= \left(\underset{1\leq s\leq n-1}{\underset{e\in B_{1,s}}\bigoplus} \mathscr{L}^{1,s}_{n-1} (e)_{(a_i,b_i)} \right)
	\bigoplus  \left(\underset{1\leq t\leq n-1}{\underset{e\in B_{2,t}}\bigoplus} \mathscr{L}^{2,t}_{n-1} (e)_{(a_i,b_i)}\right)
	\bigoplus\\ & \quad  \left(\underset{1\leq r\leq n-1}{\underset{e\in B_{3,r}}
	\bigoplus} \mathscr{L}^{1,r}_{n-1} (e,\psi(e))_{(a_i,b_i)} \right)
	\bigoplus \left(\overset{l}{\underset{j=1}\bigoplus} \mathscr{L}^{2,t}_{n-1} (h_j)_{(a_i,b_i)}\right) \bigoplus \left(\overset{m}{\underset{j=1}\bigoplus} \langle(g_j , 0),(0, g_j) \rangle\right)
\end{align*}}}
Take $B_{1,n}$ a bases for $W'_0$.
We have:
\begin{align*}
	V(a_i,b_i)&= \left(\underset{1\leq s\leq n-1}{\underset{e\in B_{1,s}}\bigoplus} \mathscr{L}^{1,s}_{n} (e)_{(a_n,b_n)} \right)
	\bigoplus  \left(\underset{1\leq t \leq n-1}{\underset{e\in B_{2,t}}\bigoplus} \mathscr{L}^{2,t}_{n} (e)_{(a_n,b_n)}\right)\bigoplus\\
	& \qquad  \left(\underset{1\leq r\leq n-1}{\underset{e\in B_{3,r}}
	\bigoplus} \mathscr{L}^{1,r}_{n} (e,\psi(e))_{(a_n,b_n)} \right)
	\bigoplus \left(\underset{e\in B_{1,n}}\bigoplus \mathscr{L}^{n}_{1,n} (e)_{(a_n,b_n)}\right) \bigoplus\\
	& \qquad\left(\overset{l}{\underset{j=1}\bigoplus} \mathscr{L}^{n}_{2,1}(h_j)_{(a_n,b_n)}\right) \bigoplus \left(\overset{m}{\underset{j=1}\bigoplus} \mathscr{L}^{n}_{3,n}(g_j, \phi(g_j))_{(a_n,b_n)}\right). 
\end{align*}
Here $\mathscr{L}^n_{3,n}(g_j , \phi(g_j))_{(a_i,b_i)} =\langle (g_j , 0),(0, g_j)\rangle$. We obtain:
\begin{align*}
	V &= \left(\underset{1\leq s\leq n}{\underset{e\in B_{1,s}}\bigoplus} \mathscr{L}^{1,s}_{n} (e)\right)\bigoplus
	\left(\underset{1\leq t\leq n-1}{\underset{e\in B_{2,t}}\bigoplus} \mathscr{L}^{2,t}_{n} (e)\right)\bigoplus \left(\underset{1\leq r\leq n-1}{\underset{e\in B_{3,r}}
	\bigoplus} \mathscr{L}^{n}_{3,r} (e,\psi(e)\right)\bigoplus\\
	&\qquad \left(\overset{l}{\underset{j=1}\bigoplus} \mathscr{L}^{n}_{2,1}(h_j)\right)\bigoplus\left(\overset{m}{\underset{j=1}\bigoplus} \mathscr{L}^{n}_{3,n}(g_j, \phi(g_j))\right)\bigoplus S(W_0).
\end{align*}

This proves our result.\\

Using the notation in \textbf{a, b} and \textbf{c} of  \ref{sect5.1}, by means of a direct calculation, the following
Lemmas can be proved.
\begin{lema}
	\begin{description}
		\item{(a.)} For all $i,\quad \textnormal{Hom}(\mathcal{L}_{3,i},\mathcal{L}_{3,i})=\left < \phi=
		\begin{pmatrix}
			1&0\\
			0& 1
		\end{pmatrix}; \quad \psi= \begin{pmatrix}
			0&1\\
			0& 0
		\end{pmatrix}
		\right >.$
		\item{(b.)} For $i > j,\quad \textnormal{Hom}(\mathcal{L}_{3,i},\mathcal{L}_{3,j})= \left < \phi=
		\begin{pmatrix}
			0&1\\
			0& 0
		\end{pmatrix};
		\psi= \begin{pmatrix}
			0&0\\
			0& 1
		\end{pmatrix}
		\right >.$
		\item{(c.)} For $i < j,\quad \textnormal{Hom}(\mathcal{L}_{3,i},\mathcal{L}_{3,j})= \left < \phi=
		\begin{pmatrix}
			1&0\\
			0& 0
		\end{pmatrix};
		\psi=
		\begin{pmatrix}
			0&1\\
			0& 0
		\end{pmatrix}
		\right >.$
	\end{description}
\end{lema}
\begin{lema}
	\begin{description}
		\item{(a.)} If $i > j,\quad  \textnormal{Hom}(\mathcal{L}_{1,i};\mathcal{L}_{3,j} )=\left < \phi =
		\begin{pmatrix}
			0\\
			1
		\end{pmatrix};
		\psi=
		\begin{pmatrix}
			1\\
			0
		\end{pmatrix}
		\right >.$
		\item{(b.)} If $i \leq j,\quad \textnormal{Hom}(\mathcal{L}_{1,i},\mathcal{L}_{3,j})=\left < \phi = 
		\begin{pmatrix}
			1\\
			0
		\end{pmatrix}
		\right >$
	\end{description}
\end{lema}

\begin{lema}
	If $\textnormal{Hom}(\mathcal{L}_{2,i},\mathcal{L}_{3,j}) \neq 0$, then $i < j$ and this $k$-vector space is generated by
	$\phi = 
	\begin{pmatrix}
		1\\
		0
	\end{pmatrix}.$
\end{lema}
\begin{teor}
	\begin{enumerate}
		\item  A representative set of the isomorphism classes of indecomposable projectives in
		$\rep(\mathfrak{U}_n)$ is given by $\mathcal{L}_{3,r}$ with $1 \leq r \leq n$ and the trivial representation $S$.
		\item  A representative set of the isomorphism classes of indecomposable injectives in $\rep(\mathfrak{U}_n)$ are given by $\mathcal{L}_{3,r}$ with  $1 \leq r \leq n$ and $\mathcal{L}_{2,n}$.
		\item  The following are all the almost split sequences in $\rep(\mathfrak{U}_n)$:
		
		$$\mathscr{L}_{1,1}\xrightarrow{(0,1)^t} \mathscr{L}_{3,1}\xrightarrow{(1,0)}  \mathscr{L}_{1,1}.$$
		For $2 \leq  i \leq  n$, $$\mathscr{L}_{1,i}\xrightarrow{(0,1)^t} \mathscr{L}_{3,i-1}\xrightarrow{(1,0)} \mathscr{L}_{2,i}$$
		For $1 \leq i \leq n-1$,  $$\mathscr{L}_{2,i}\xrightarrow{(0,1)^t} \mathscr{L}_{3,i+1}\xrightarrow{(1,0)} \mathscr{L}_{1,i}$$
		and  $$ {S}\xrightarrow{(0,1)^t} \mathscr{L}_{3,n}\xrightarrow{(1,0)} \mathscr{L}_{2,n}.$$
	\end{enumerate}
\end{teor}
\textbf{Proof}
\begin{enumerate}
	\item  By Proposition 5, $\mathscr{L}_{3,i} = P(\{a_i
	, b_i\})$, and $S$ are all the projectives.
	\item  Observe that ${\mathfrak{U}}_n^{\textnormal{op}}$
	is also of the form ${\mathfrak{U}}_n$, then indecomposables $W\in  \rep({\mathfrak{U}}_n^{\textnormal{op}})$ with
	$\textnormal{dim}_k W_0 = 2$ are projectives, this implies by Corollary 3, that indecomposables $V$ in $\rep(\mathfrak{U}_n)$ with $\textnormal{dim}_kV_0 = 2$ are injectives, therefore all the $\mathscr{L}_{3,i}$  are injectives, and
	$\mathscr{L}_{2,n}$  coincides with the representation $J$ of Remark 4, therefore the injectives are the
	projectives $\mathscr{L}_{3,i}$ and $\mathscr{L}_{2,n}$. 
	\item All sequences in this item are $\varepsilon$-sequences. We will prove that the morphism
	$\mathscr{L}_{1,i}
	\xrightarrow{(0,1)^t} \mathscr{L}_{3,i-1}$ is a left almost split morphism. For this it is enough to prove that for any
	morphism $v:\mathscr{L}_{1,i}\to V$, with $V$ indecomposable non isomorphic to $\mathscr{L}_{1,i}$ there is a
	morphism $u: \mathscr{L}_{3,i-1}\to V$ with $u(0,1)^t = v$. We may assume $v\neq 0$, and $v(e) = e$.
	Then if $V = \mathscr{L}_{1,j}$ we must have $j< i$, and the linear map $u=(0, 1): k \langle e_1, e_2 \rangle \to k \langle e \rangle$ is a morphism from $\mathscr{L}_{3,i-1}\to  \mathscr{L}_{1,j}$ with $u(0, 1)^t = v$. Now if $V = \mathscr{L}_{2,j}$, again
	we may suppose that $v(e) = e$. Then the linear map $u = (0, 1): k \langle e_1, e_2 \rangle \to ke$
	gives a morphism from $\mathscr{L}_{3,i-1}\to \mathscr{L}_{1,i}$ with $u(0, 1)^t = v$. Now suppose $V = \mathscr{L}_{3,j}$ , if
	$i-1 = j$, then $v = (d, c)^t$. The linear map $u =
	\begin{pmatrix}
		c& d\\
		0&c
	\end{pmatrix}: k \langle e_1, e_2 \rangle \to k \langle e_1, e_2 \rangle$
	is a morphism from $\mathscr{L}_{3,i-1}\to \mathscr{L}_{3,i-1}$ such that $u(0, 1)^t = v$. Suppose now $i-1 > j$, the linear map $u=\begin{pmatrix}
		0& c\\
		0&d
	\end{pmatrix}$
	is a morphism $\mathscr{L}_{3,i-1} \to  \mathscr{L}_{3,j}$ with $u(0, 1)^t = v$. In
	case $i < j$, we have $v = (c, 0)^t$, and there is a morphism $u: \mathscr{L}_{3,i-1}\to  \mathscr{L}_{3,j}$ given by
	the linear map
	$u=\begin{pmatrix}
		0& c\\
		0&0
	\end{pmatrix}$, such that $v = u(0, 1)^t$. In a similar way, using the
	above lemmas one can prove that the first morphism in each of the sequences is a
	left almost split morphism. Since the other end of the sequence is indecomposable,
	therefore all these sequences are almost split sequences. 
\end{enumerate}
In this way we obtain that the Auslander-Reiten quiver for $\textnormal{D}_n$ is given by
{\tiny{
		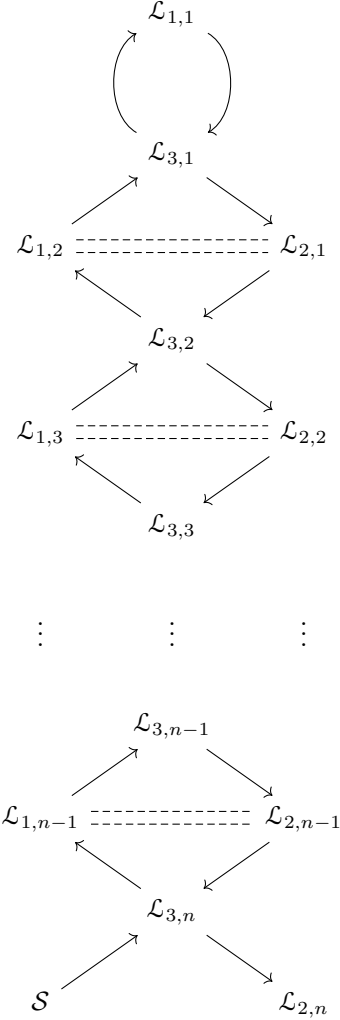
\begin{figure}[H]
			\begin{center}
				\begin{tikzcd}[column sep= small,row sep=normal]
					&\mathscr{L}_{1,1}  \arrow[dd,bend left =60,looseness=0.9]  &\\
					&&&	\\
					&\mathscr{L}_{3,1} \arrow[uu,bend left =60,looseness=0.9]\arrow[dr]&\\ \mathscr{L}_{1,2}\arrow[ur]\arrow[dash, dashed, shift right, rr] \arrow[dash, dashed, shift left, rr]&  &\mathscr{L}_{2,1}\arrow[dl]\\
					&\mathscr{L}_{3,2}\arrow[ul] \arrow[dr]&\\
					\mathscr{L}_{1,3}\arrow[ur]\arrow[dash, dashed, shift right, rr] \arrow[dash, dashed, shift left, rr]&  &\mathscr{L}_{2,2}\arrow[dl]\\
					&\mathscr{L}_{3,3} \arrow[ul]&\\
					\vdots& \vdots& \vdots\\
					&\mathscr{L}_{3,n-1}\arrow[dr] &\\
					\mathscr{L}_{1,n-1}\arrow[ur]\arrow[dash, dashed, shift right, rr] \arrow[dash, dashed, shift left, rr] &  &\mathscr{L}_{2,n-1}\arrow[dl]\\
					&\mathscr{L}_{3,n} \arrow[ul]\arrow[dr]&\\
					\mathcal{S} \arrow[ur]&  &\mathscr{L}_{2,n}\\
				\end{tikzcd}
				
				\caption{Auslander-Reiten quiver of poset with involution of type $\text{D}_n$.}
			\end{center}
\end{figure}}}
\subsection{Poset with an involution of type $\mathfrak{U}_\infty$}
The previous results can be extended to posets of type $\mathfrak{U}_\infty$ in the following way. We define the functor
$$\Xi_n : \rep(\mathfrak{U}_n) \to  \rep(\mathfrak{U}_\infty)$$ such that 
if  $V = (V_0, V_{(a_i,b_i)})_{1\leq i\leq n}$  then $\Xi(V) = (V_0, V_{(a_i,b_i)})_{1\leq i}$ with  $V_{(a_j,b_j )} =
(0, V_0)$  for  $j > n$.
\par\bigskip
Clearly  if  $f: V \to W$ is a morphism in $ \rep(\mathfrak{U}_n)$ determined by the morphism 
$f: V_0 \to W_0$ then this morphism  is the same in the category $\rep(\mathfrak{U}_\infty)$.  So,  taking $\Xi _n(f) = f$  we have defined the functor $\Xi_n$.
Analogously we define a functor:
$\Theta_{n}:\rep(\mathfrak{U}_{n})\to \rep(\mathfrak{U}_{n+1})$ and it is obtained that $\Xi_{n+1}\Theta_n = \Xi_n$.
\smallskip
\begin{prop}
	The indecomposable representations of $\mathfrak{U}_{\infty}$ are the representations 
	
	$\hat{L}_{s,i} = \Xi_i(L_{s,i})$ together with
	the  simple trivial representation $S$.
\end{prop}
\textbf{Proof.} Let  $V =(V_0, V_{(a_i,b_i)})_{1\leq i}$ be an indecomposable representation of  $\mathfrak{U}_\infty$. We suppose that $V_{(a_1,b_1)} = 0$, then for each  $n$ the restriction of $V$ to $\mathfrak{U}_n$ is the form $(V)_{\mathfrak{U}_n} = \underset{i=1}{\overset{s}\bigoplus}{L}_s$ with  ${L}_s = \mathscr{L}_{1,i}$ with $i > 1$ or the trivial representation. Since $V_0$ is finite dimensional there exists   $n$  and finite sum $W = \underset{i=1}{\overset{s}\bigoplus}{L}_s$
in $\mathfrak{U}_n$ such that for all $m > n$,  $V$ restricted to $\mathfrak{U}_m$ coincides with the restriction of  $\Xi_n(W)$ to $\mathfrak{U}_m$, therefore $V = \Xi_n(W)$  and since $V$ is an   indecomposable, then $V = \Xi_n(\mathscr{L}_{1,n})$.\\

Now we suppose that $V_{(a_1,b_1)}\neq 0$ then $V$ restricted to $\mathfrak{U}_n$ is the form $\underset{i=1}{\overset{s}\bigoplus}{L}_s$ where each $L_s$ is the form $\mathscr{L}_{j,i}$ with $j = 1, 2, 3$ and  at least one $L_s$ has the form $\mathscr{L}_{2,i}$ or $\mathscr{L}_{3,i}$. As before, there exists $n$ such that for all $m \geq n$ the restriction from $V$ to
$\mathfrak{U}_m$ coincides with the restriction of $\Xi_n(W)$ to $\mathfrak{U}_m$ therefore $V=\Xi_n(\mathscr{L}_{j,n})$ with
$j = 2$ or $j=3$.
$\hfill\blacksquare$
\smallskip
\begin{prop}
	Let $a: X \xrightarrow{u} Y \xrightarrow{v} Z$
	be an almost split sequence in 
	$\rep(\mathfrak{U}_n)$ with $X$ different from the  trivial representation. Then 
	$$b:\Xi_n(X) \xrightarrow{\Xi_n(u)} \Xi_n(Y) \xrightarrow{\Xi_n(v)} \Xi_n(Z),$$
	is an almost split sequence in  $\rep(\mathfrak{U}_\infty)$.
\end{prop}

\textbf{Proof.}  The sequence  $b$ is a nontrivial  $\varepsilon$-sequence in  $\rep(\mathfrak{U}_\infty)$  whose 
extremes are indecomposable. Let $h: Y\to \Xi_n(Y)$ be a morphism that is not  a retraction in $\rep(\mathfrak{U}_\infty)$ with $Y$ indecomposable, then $Y = \Xi_m(W)$ with $W$ indecomposable in
$\rep(\mathfrak{U}_m)$ for some $m > n$. We have that $\Xi_m(a)$ is an almost split  sequence in  $\rep(\mathfrak{U}_m)$, so we can suppose that $m = n$ and then $h = \Xi_m(w)$ where  $w: W \to Z$ is  a morphism that is not a  retraction.  Therefore,  there exists $g: W \to Y$ with  $vg = w$, thus $h = \Xi_m(v)\Xi_m(g)$. This proves our assertion. 

In this way, we obtain that the Auslander-Reiten quiver for $\text{D}_\infty$ is given by
\enlargethispage{5cm}
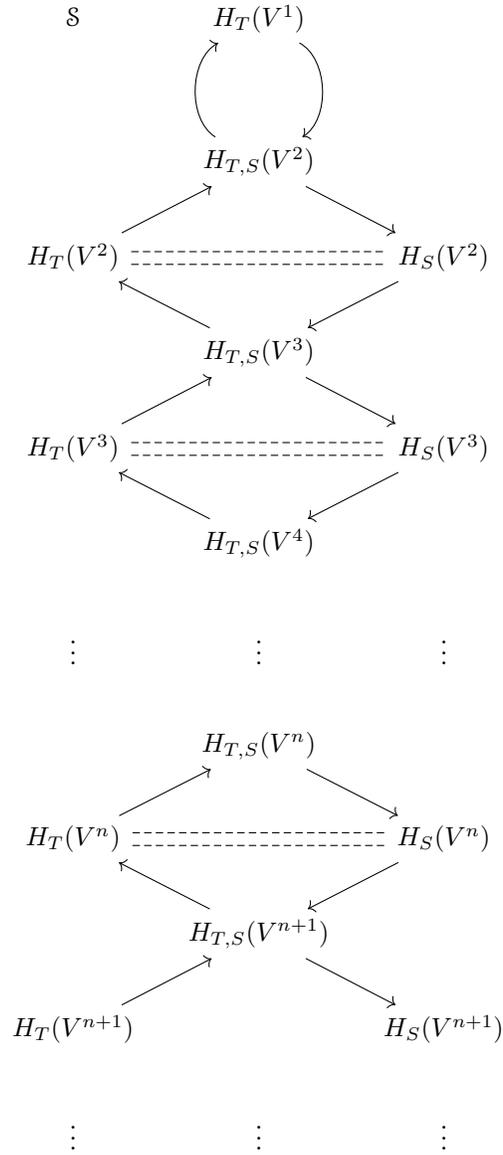
\begin{figure}[H]
	\begin{center}
		\begin{tikzcd}[column sep= small,row sep=normal]
			\mathscr{S}	&H_T(V^1)  \arrow[dd,bend left =60,looseness=0.9]  &\\
			&&&	\\
			&H_{T,S}(V^2) \arrow[uu,bend left =60,looseness=0.9]\arrow[dr]&\\ H_{T}(V^2)\arrow[ur]\arrow[dash, dashed, shift right, rr] \arrow[dash, dashed, shift left, rr]&  &H_{S}(V^2)\arrow[dl]\\
			&H_{T,S}(V^3)\arrow[ul] \arrow[dr]&\\
			H_{T}(V^3)\arrow[ur]\arrow[dash, dashed, shift right, rr] \arrow[dash, dashed, shift left, rr]&  &H_{S}(V^3)\arrow[dl]\\
			&H_{T,S}(V^4) \arrow[ul]&\\
			\vdots& \vdots& \vdots\\
			&H_{T,S}(V^n)\arrow[dr] &\\
			H_{T}(V^n)\arrow[ur]\arrow[dash, dashed, shift right, rr] \arrow[dash, dashed, shift left, rr] &  &H_{S}(V^n)\arrow[dl]\\
			&H_{T,S}(V^{n+1}) \arrow[ul]\arrow[dr]&\\
			H_{T}(V^{n+1})\arrow[ur] &  &H_{S}(V^{n+1})\\
			\vdots& \vdots& \vdots\\
		\end{tikzcd}
	\end{center}
	\caption{Auslander-Reiten quiver of a poset with an  involution of type $\text{D}_\infty$.}
\end{figure}

\par\bigskip

Raymundo Bautista Ramos\\
raymundo@matmor.unam.mx \\
Universidad Nacional Autónoma de México\\
\par
Verónica Cifuentes Vargas\\
vcifuentesv@udistrital.edu.co\\
Universidad Distrital Francisco José de Caldas.

\end{document}